\documentclass[titlepagea4paper,11pt]{article}
\usepackage[arrow,matrix]{xy}
\usepackage{amssymb}
\usepackage{amsbsy,amsmath}
\newtheorem{thm}{Theorem}[section]
\newtheorem{lemma}[thm]{Lemma}
\newtheorem{cor}[thm]{Corollary}

\newcommand{\ratmap}{\cdots\!\rightarrow}
\newcommand{\seq}{\longrightarrow}

\newcommand{\bb}{\mathbb}

\newenvironment{proof}{\vspace{3pt}\indent
                       \textsc{\bf Proof:}\quad }
                       {\hfill$\square$\vspace{3pt}}
\newtheorem{example}[thm]{Example}
\newtheorem{construction}[thm]{Construction}
\newenvironment{ex}[1]{\begin{example}\rm{\bf #1}}{\hfill$\bigtriangleup$\end{example}}
\newenvironment{constr}{\begin{construction}\rm}{\hfill$\bigtriangleup$\end{construction}}

\title{Calabi-Yau-threefolds with Picard number $\rho(X)=2$ and their K\"ahler cone II}
\author{Marco K\"uhnel\footnote{The author acknowledges gratefully support by the DFG priority program 'Global Methods in Complex Geometry'.}\\  Department of Mathematics\\  University of Bayreuth\\ D-95440 Bayreuth}
\date{\today}

\begin{document}

\maketitle

\begin{abstract}We prove the rationality of the K\"ahler cone and the positivity of $c_2(X)$, if
$X$ is a Calabi-Yau-threefold with $\rho(X)=2$ and has an embedding into a ${\bb P}^n$-bundle over
${\bb P}^m$ in the cases $(n,m)=(1,3),(3,1)$. The case $(n,m)=(2,2)$ has been done in the first part
of this paper. Moreover, if $(n,m)=(3,1)$, we describe the 'other' contraction different from the 
projection.
\end{abstract}

\section{Introduction}In this paper, a Calabi-Yau-threefold is a compact complex K\"ahler manifold of dimension three with $K_X={\cal O}_X$ and
$H^1({\cal O}_X)=0$. 

Wilson stated in 1994 \cite{wi94} a conjecture about the rationality of the K\"ahler cone of a Calabi-Yau-threefold. It says that the K\"ahler cone
of a Calabi-Yau-threefold $X$ is rational and finitely generated in $N^1(X)$, if $c_2(X)$ is positive, i.e. $D.c_2(X)>0$ for every nef divisor $D$.

In \cite{k01} we dealt with the case $\rho(X)=2$. First we proved some general results about the
K\"ahler cone and then
concentrated on the case that $X$ was embedded in a ${\bb P}^2$-bundle over 
${\bb P}^2$. For this class of Calabi-Yau-manifolds we confirmed Wilson's conjecture. 

In this paper, we want to finish this track by considering $X$ which are embedded in
either a ${\bb P}^1$-bundle over ${\bb P}^3$ or a ${\bb P}^3$-bundle over ${\bb P}^1$.
The first case offers some interesting perspectives.
The Calabi-Yau-manifolds turn out to be generic double covers of ${\bb P}^3$ ramified over
an octic. We compute the number of fibres of the bundle projection in $X$ and describe the K\"ahler cone. The 
ramifying octics deserve a more detailed discussion in another paper.
  
Since $C={\bb P}^1$, if $X\seq C$ is a fibration onto a normal curve $C$ and $X$ a Calabi-Yau manifold,
it is also natural to turn our attention to those $X$, which can be embedded an ${\bb P}^3$-bundle over
${\bb P}^1$. 

We will denote $K(X)$ for the K\"ahler cone of $X$. 
The main results are

\begin{thm}
Let $X\subset{\bb P}({\cal E})$ be a Calabi-Yau-$3$-fold with $\rho(X)=2$, 
with ${\cal E}$ either being a rank-$2$-bundle over ${\bb P}^3$ such that
$h^0(-K_{{\bb P}({\cal E})})>1$ or an arbitrary  rank-$4$-bundle over ${\bb P}^1$. Then
$\partial K(X)$ is rational and semiample. Furthermore all $D\in\overline{K(X)}$ satisfy
$D.c_2(X)>0$.
\end{thm}

In the case of ${\cal E}$ being a rank-4-bundle over ${\bb P}^1$ we can give a complete classification 
of the occuring contractions.

\begin{thm}Let $X\subset{\bb P}({\cal E}):=Z$ be a Calabi-Yau-3-fold with $\rho(X)=2$, where
${\cal E}\seq{\bb P}^1$ is of rank 4. Let us fix 
${\cal E}={\cal O}\oplus{\cal O}(a_1)\oplus{\cal O}(a_2)\oplus{\cal O}(a_3)$, with $0\le a_1\le a_2\le a_3$
and let $\psi:X\seq X'$ be the second contraction. If $E$ denotes the exceptional locus of $\psi$ and 
${\cal F}:={\cal O}\bigoplus_{i|a_i=0}{\cal O}(a_i)$ the maximal trivial subbundle of ${\cal E}$, then: 
\begin{enumerate}
\item If $c_1({\cal E})=3$, then $rk{\cal F}\le 2$ and $E={\bb P}({\cal F})\cong{\bb P}^1\times{\bb P}^{rk{\cal F}-1}$.
\item If $c_1({\cal E})=2$, then $rk{\cal F}\in\{2,3\}$ und $E=X\cap{\bb P}({\cal F})={\bb P}^1\times Y$, with $\dim Y=rk{\cal F}-2.$ 
\begin{enumerate}
\item If $rk{\cal F}=2$, then $Y$ consists of four points.
\item If $rk{\cal F}=3$, then $Y$ is a smooth plane quartic.
\end{enumerate}
\item If $c_1({\cal E})=1$, then $Z$ is the blow-up of ${\bb P}^4$ in a linearly embedded ${\bb P}^2$; if $X\in|-K_Z|$ is general, then
$E=\bigcup_{i=1}^{16} C_i$, with $C_i\cong{\bb P}^1$; furthermore, $X'$ is a quintic in ${\bb P}^4$ with
 $16$ double points on a linearly embedded
${\bb P}^2$.
\item If $c_1({\cal E})=0$, then $Z={\bb P}^1\times{\bb P}^3$ and $E=\bigcup_{i=1}^{64}C_i,$ with
 $C_i\cong{\bb P}^1$.
\end{enumerate}
\end{thm}

For the proofs of these results we proceed like in \cite{k01} and prove first a generalization of a
lemma of Koll\'ar \cite{bo89}:

\begin{thm}Let $X\subset{\bb P}({\cal E})$ be a Calabi-Yau-threefold with $\rho(X)=2$, 
with ${\cal E}$ either being a rank-$2$-bundle over ${\bb P}^3$ or a rank-$4$-bundle over ${\bb P}^1$.
Then $K(X)=K(Z)|X$.
\end{thm}

In contrast to the generalized Koll\'ar Lemma in \cite{k01}, the situation here is simpler, since
there are no exceptions in the theorem above. However, in the case of ${\cal E}$ being a rank-2-bundle
over ${\bb P}^3$ it is has to be investigated, if $\pi^*h$ is not ample. It turns out,
that this is true, if $\rho(X)>1$. So the proofs of the rationality
result and the positivity of $c_2(X)$ get much shorter than those in \cite{k01}.
In the case of ${\cal E}$ being a rank-4-bundle over ${\bb P}^1$, the rationality of the K\"ahler cone
can even be proven without using the generalized Koll\'ar Lemma.

This article grew out of the author's doctoral thesis at the University of Bayreuth.

\section{Notation}

In this section we summarize the most important notations of this paper. $X$ will always denote a Calabi-Yau-threefold, while $Z$ is
a fourfold in every case. 
\begin{center}
\begin{tabular}{|l|l|}
\hline
$N^1(X)$& the ${\bb R}$-vector space of numerical classes of $Div(X)\otimes{\bb R}$\\
$K(X)$& the open part of the K\"ahler cone of $X$, i.e. the ample cone\\
$W(X)$& the hypersurface $\{D^3=0\}\subset N^1(X)$\\
${\cal O}_X(1)$ & the restriction ${\cal O}_Z(1)|X$, where ${\cal O}_Z(1)$ is the tautological bundle\\
& associated to $Z={\bb P}({\cal E})$\\
$K_Z$& the canonical divisor of $Z$\\
$c_i({\cal E})$& the $i$-th Chern class of a bundle ${\cal E}$\\
$c_i(M)$& the $i$-th Chern class of the tangent bundle of the complex\\
& manifold $M$\\
\hline
\end{tabular}
\end{center}

\section{Some general statements}

For a more detailed description of properties of the K\"ahler cone of a Calabi-Yau-threefold
with $\rho(X)=2$ we refer to \cite{diss,k01}. At this place, we mention only results which are necessary
for the later parts of this article.

The first fundamental theorem we want to cite is proved by Wilson in \cite{wi}.

\begin{thm}{\bf (Wilson)} Let $X$ be a Calabi-Yau-threefold. If $D\in\partial K(X)$ and $D^3>0$, then
there is some $r\in{\bb R}$ with $rD\in Pic(X)$.
\end{thm}

Hence it is natural to condider the cubic hypersurface $W(X):=\{D\in N^1(X)|D^3=0\}$. A useful
statement is

\begin{lemma}Let $X$ be a Calabi-Yau-threefold with $\rho(X)=2$. If $W(X)$ contains a double line,
then $W(X)$ is rational.
\end{lemma}

\begin{proof}$W(X)$ is an appropriate affine neighbourhood of ${\bb P}(N^1(X))\cong{\bb P}^1$
given by some cubic polynomial $w\in{\bb Z}[x]$. Let denote
$Dw\in{\bb Q}[x]$ the formal derivative of $w$. Then, if
$$w=(x-a)^2(x-b)$$
and $a\not=b$
$$(x-a)=gcd(w,Dw)\in{\bb Q}[x].$$
if $a=b$, then 
$$(x-a)^2=gcd(w,Dw)\in{\bb Q}[x].$$
In both cases, $a,b\in{\bb Q}$ follows.
\end{proof}

Putting both results together, we get in particular

\begin{cor}\label{p1}If $X$ is a Calabi-Yau-manifold with $\rho(X)=2$ and $\phi:X\seq{\bb P}^1$
a fibration, then $\partial K(X)$ is rational.
\end{cor}

Finally, we need a Lemma of Koll\'ar, proved in \cite{bo89}.

\begin{lemma}\label{0kollar}{\bf (Koll\'ar)}
If $Z$ is a Fano-4-fold and $X\in|-K_Z|$ a Calabi-Yau-manifold, then 
$i^*:N^1(Z)\seq N^1(X)$
is an isomorphism and $K(X)=K(Z)|X$.
\end{lemma}

\section{Calabi-Yau-threefolds in ${\bb P}^1$-bundles over ${\bb P}^3$}

We are interested in Calabi-Yau-threefolds $X$ of the form
$X\subset{\bb P}({\cal E})=:Z$, $X\in|-K_{{\bb P}({\cal E})}|$, with ${\cal E}$ being a vectorbundle of 
rank $2$ over 
${\bb P}^3$. Let us denote $p:{\bb P}({\cal E})\seq {\bb P}^3$ the bundle projection and 
$\pi:X\seq{\bb P}^3$ the restriction of $p$ to $X$. 
The hyperplane class in ${\bb P}^2$ shall be denoted by $h$, the fibre of $p$ by $F$. The expression
$\gamma({\cal E}):=(c_1^2({\cal E})-4c_2({\cal E})).h$ is invariant under ${\cal E}\mapsto{\cal E}\otimes L$, with $L$ being a line bundle
over ${\bb P}^2$. The line bundle ${\cal O}_Z(1)|X$ will be called ${\cal O}_X(1)$.

The following sequences are basic for our proofs and results:

\begin{eqnarray}
&0\seq T_{Z|{\bb P}^3}\seq T_Z \seq p^* T_{{\bb P}^3}\seq 0&\label{e1}\\
&0\seq {\cal O}_Z\seq p^*({\cal E}^\vee)\otimes{\cal O}_Z(1) \seq T_{Z|{\bb P}^3}\seq 0&\label{e2}\\
&0\seq T_X\seq T_Z|X \seq N_{X|Z}\seq 0&\label{e3}\\
&{\cal O}_Z(1)^2-p^*c_1({\cal E}).{\cal O}_z(1)+p^*c_2({\cal E})=0&\label{e4}
\end{eqnarray}


 

By the K\"unneth formula we get $b_1(Z)=0, b_2(Z)=2, b_3(Z)=0, b_4(Z)=2$.

The intersection theory on $Z$ is computed inductively by ${\cal O}_Z(1).p^*h^3=1$ and equation
(\ref{e4}):

\begin{lemma}\label{3intersection}Let ${\cal E}\seq{\bb P}^3$ be a rank-2-bundle and 
$Z:={\bb P}({\cal E})$. Then holds:
\begin{enumerate}
\item ${\cal O}_Z(1).p^*h^3=1$,
\item ${\cal O}_Z(1)^2.p^*h^2=c_1({\cal E}).h^2$,
\item ${\cal O}_Z(1)^3.p^*h=c_1^2({\cal E}).h-c_2({\cal E})$,
\item ${\cal O}_Z(1)^4=c_1^3({\cal E})-2c_1({\cal E})c_2({\cal E}).$
\end{enumerate}
\end{lemma}

\subsection{Intersection Product and Picard Number}

By standard computations we get

\begin{lemma}\label{3c2}
Sei $X\subset{\bb P}({\cal E})$ eine Calabi-Yau-Mannig\-faltig\-keit, 
wobei ${\cal E}\seq{\bb P}^3$ ein Rang-2-Vektorb\"undel ist. Dann gilt
\begin{enumerate}
\item $c_3(X)=-8\gamma -168$
\item $\pi^*h.c_2(X)=44$
\item ${\cal O}_Z(1)|X.c_2(X)=4\gamma+22c_1({\cal E})+24$
\item $-K_Z|X.c_2(X)=8\gamma +224$
\item $\pi^*h^3=2$
\item ${\cal O}_Z(1)|X.\pi^*h^2=c_1({\cal E}).h^2+4$
\item ${\cal O}_Z(1)^2|X.\pi^*h=\frac 12\gamma +\frac 12c_1^2({\cal E}).h+4c_1({\cal E}).h^2$
\item ${\cal O}_Z(1)^3|X=\gamma+\frac 34\gamma c_1({\cal E}).h^2+3c_1^2({\cal E}).h+\frac 14c_1^3({\cal E})$
\end{enumerate}
\end{lemma}

Compare the following result about the Picard number to the corresponding theorem in \cite{k01}.

\begin{thm}\label{3rho}Let $X\subset{\bb P}({\cal E})$ be a Calabi-Yau-manifold, 
with ${\cal E}\seq{\bb P}^3$ being a rank-2-bundle. 
If ${\cal E}$ is stable and $H^1(-K_Z)=H^2(-K_Z)=0$ (e.g. $-K_Z$ big and nef), then
$$\rho(X)=2+h^2({\cal E}^\vee\otimes{\cal E}).$$
\end{thm}

\begin{proof}
We look at the sequences
\begin{eqnarray}&0\seq N_{X|Z}^\vee\seq\Omega_Z|X\seq\Omega_X\seq 0&\label{3?}\end{eqnarray}
and
\begin{eqnarray}&0\seq\Omega_Z\otimes K_Z\seq\Omega_Z\seq\Omega_Z|X\seq 0.&\label{3??}\end{eqnarray}

First, we want to show $H^i(T_Z)=H^i({\cal E}^\vee\otimes{\cal E})$ for $i>1$. For this purpose we
compute
$R^ip_*(p^*({\cal E}^\vee)\otimes{\cal O}_Z(1))={\cal E}^\vee\otimes R^ip_*{\cal O}_Z(1)=0$ for $i>0$. Hence
the Leray spectral sequence implies 
$$H^i(p^*({\cal E}^\vee)\otimes{\cal O}_Z(1))=H^i({\cal E}^\vee\otimes{\cal E}).$$
Sequence (\ref{e2}) shows, that $H^i(p^*({\cal E}^\vee)\otimes{\cal O}_Z(1))=H^i(T_{Z|{\bb P}^3})$ for $i>0$,
since $R^ip_*{\cal O}_Z=0$ for $i>0$ and therefore 
$H^i({\cal O}_Z)=H^i({\cal O}_{{\bb P}^3})=0$ for $i>0$.
To apply sequence (\ref{e1}), we compute
$R^ip_*p^*T_{{\bb P}^3}=T_{{\bb P}^3}\otimes R^ip_*{\cal O}_Z=0$ for $i>0$. 
Therefore we see again by the Leray spectral sequence 
$$H^i(p^*T_{{\bb P}^3})=H^i(T_{{\bb P}^3})=0$$
for $i>0$. This implies by sequence (\ref{e1}), that
$$H^i(T_{Z|{\bb P}^3})=H^i(T_Z)$$
for $i>1$. Hence
$$H^i(T_Z)=H^i({\cal E}^\vee\otimes{\cal E})$$
for $i>1$.

We now know 
$$H^i(\Omega_Z\otimes K_Z)=H^{4-i}(T_Z)^\vee=H^{4-i}({\cal E}^\vee\otimes{\cal E})^\vee$$
for $i<3$. In particular,
$$H^1(\Omega_Z\otimes K_Z)=H^3({\cal E}^\vee\otimes{\cal E})^\vee=H^0({\cal E}^\vee\otimes{\cal E}\otimes{\cal O}(-4))=0,$$
with the last equality holding since ${\cal E}$ is stable and therefore simple.
From $b_3(Z)=0$ follows $H^{2,1}(Z)=0$. Since $N_{X|Z}^\vee=K_Z|X$ the cohomology sequences of (\ref{3??})
and (\ref{3?}) contain  
\begin{eqnarray}&0\seq H^1(\Omega_Z)\seq H^1(\Omega_Z|X)\seq H^2({\cal E}^\vee\otimes{\cal E})^\vee\seq 0&\label{3!}
\end{eqnarray}
resp.
\begin{eqnarray}&0\seq H^1(K_Z|X)\seq H^1(\Omega_Z|X)\seq H^1(\Omega_X)\seq H^2(K_Z|X).&\label{3!!}\end{eqnarray}

By the assumption $H^1(-K_Z)=H^2(-K_Z)=0$ we conclude
by
$$0\seq 2K_Z\seq K_Z\seq K_Z|X\seq 0$$
and 
\begin{eqnarray}
H^1(K_Z)=&H^3({\cal O}_Z)^\vee=&0\nonumber\\
H^2(K_Z)=&H^2({\cal O}_Z)^\vee=&0\nonumber\\
H^2(2K_Z)=&H^2(-K_Z)^\vee=&0\nonumber\\
H^3(2K_Z)=&H^1(-K_Z)^\vee=&0,\nonumber
\end{eqnarray}
that 
$$H^1(K_Z|X)=H^2(K_Z|X)=0.$$

Therefore by (\ref{3!!}) and (\ref{3!})
$$\rho(X)=h^1(\Omega_Z|X)=\rho(Z)+h^2({\cal E}^\vee\otimes{\cal E}).$$
\end{proof}

\begin{ex}{(A Calabi-Yau manifold with $\rho(X)=1$)}\label{3rho=1}
Let ${\cal E}={\cal O}\oplus{\cal O}(4)$. Then we prove like in Thm \ref{3rho}, that 
$$H^i(\Omega_Z\otimes K_Z)=H^{4-i}(T_Z)^\vee=H^{4-i}({\cal E}^\vee\otimes{\cal E})^\vee.$$
But now
$$H^3({\cal E}^\vee\otimes{\cal E})^\vee=H^0({\cal E}^\vee\otimes{\cal E}\otimes{\cal O}(-4))=H^0({\cal O}(-8)\oplus
2{\cal O}(-4)\oplus{\cal O})={\bb C}.$$

The cohomology sequence of (\ref{3?}) starts
$$0\seq H^0(N_{X|Z}^\vee)\seq H^0(\Omega_Z|X)\seq H^0(\Omega_X).$$
Since $H^0(\Omega_X)=0$ and $H^0(N_{X|Z}^\vee)=0$, we see
$$H^0(\Omega_Z|X)=0.$$
This we use in the cohomology sequence of (\ref{3??}) and get the sequence
$$0\seq{\bb C}\seq H^1(\Omega_Z)\seq H^1(\Omega_Z|X)\seq H^2({\cal E}^\vee\otimes{\cal E})^\vee\seq 0.$$
We know $H^2({\cal E}^\vee\otimes{\cal E})^\vee=0$. Since $-K_Z={\cal O}_Z(2)$ is big and nef, 
$H^i(N_{X|Z}^\vee)=H^i(K_Z|X)=0$
follows for $i=1,2$. By using the cohomology sequence of (\ref{3?}) we get
$$H^1(\Omega_X)=H^1(\Omega_Z|X).$$ 
This finally implies
$$\rho(X)=1.$$
\end{ex}

As a last subject in this section, we are interested in some bounds for $\gamma$. This yields
a total lower bound for $c_3(X)$ of the here considered Calabi-Yau-threefolds, and, what is more important
within this framework, allows us to compute the number of full fibres of $p$ contained in $X$. 
The latter will be done in the next section. 

\begin{lemma}\label{3maxg}Let ${\cal E}\seq{\bb P}^3$ be a rank-2-bundle. Denote its generic
splitting type by $(a,b)$. Then
$$\gamma({\cal E})\le\gamma({\cal O}(a)\oplus{\cal O}(b))=(a-b)^2.$$
Moreover, equality holds if and only if ${\cal E}={\cal O}(a)\oplus{\cal O}(b)$.
\end{lemma}

\begin{proof}
Let us denote ${\cal E}':={\cal O}(a)\oplus{\cal O}(b)$.
By tensorizing with ${\cal O}(m)$ for $m\gg 0$, we may assume that 
$H^i({\cal E})=0$ for $i>0$. For a general hyperplane we
$H\subset{\bb P}^3$ we look at the sequence
$$0\seq{\cal E}(-1)\seq{\cal E}\seq{\cal E}|H\seq 0.$$
We see that
$$h^0({\cal E})\le h^0({\cal E}(-1))+h^0({\cal E}|H).$$
The same argument shows inductively
$$h^0({\cal E})\le h^0({\cal E}(-k))+\sum_{i=0}^{k-1}h^0({\cal E}(-i)|H)$$
for all $k$ and hence
$$h^0({\cal E})\le\sum_{i\ge 0}h^0({\cal E}(-i)|H).$$
We choose a general line $L\subset H$, and conclude by replacing ${\bb P}^3$ by $H$,
that
$$h^0({\cal E}(-i)|H)\le\sum_{j\ge 0}h^0({\cal E}(-i-j)|L).$$
Therefore 
\begin{eqnarray}h^0({\cal E})&\le&\sum_{0\le i,j}h^0({\cal E}(-i-j)|L)\nonumber\\
&=&\sum_{0\le j\le a,i+j\le a}a+1-i-j+\sum_{0\le j\le b,i+j\le b}b+1-i-j.\nonumber\\
&=&\sum_{j=0}^a{a+2-j \choose 2}+\sum_{j=0}^b{b+2-j \choose 2}\nonumber\\
&=&\sum_{j=2}^{a+2}{j\choose 2}+\sum_{j=2}^{b+2}{j\choose 2}\nonumber\\
&=&{a+3\choose 3}+{b+3\choose 3}\nonumber\\
&=& h^0({\cal E}').\nonumber
\end{eqnarray}

By assumption we have
$$\chi({\cal E})=h^0({\cal E})\le h^0({\cal E}')=\chi({\cal E}').$$
Using Riemann-Roch, thie inequality transforms to
\begin{eqnarray}\frac 18\gamma(c_1({\cal E}).h^2+4)+\frac{1}{24}c_1^3({\cal E})+\frac 12c_1^2({\cal E}).h+\frac{11}{6}c_1({\cal E}).h^2+2
&\le&\nonumber\\
\frac 18\gamma(c_1({\cal E}').h^2+4)+\frac{1}{24}c_1^3({\cal E}')+\frac 12c_1^2({\cal E}').h+\frac{11}{6}c_1({\cal E}').h^2+2.& &\nonumber
\end{eqnarray}
Since by assumption $c_1({\cal E})=c_1({\cal E}')$ and $c_1({\cal E}).h^2>0$, this implies
$$\gamma({\cal E})\le \gamma({\cal E}').$$

Equality holds if and only if all connecting homomorphisms
$$H^0({\cal E}(-i)|L)\seq H^1({\cal E}(-i-1)|H)$$
are the zero map, what means that
$$h^1({\cal E}(-i-1)|H)\le h^1({\cal E}(-i)|H)$$
for all $i\ge 0$. This implies
$$H^1({\cal E}(n)|H)=0 \mbox{ for all } n\in{\bb Z}$$
by the choice of $m$.

We conclude with Horrocks' splitting criterion, that ${\cal E}|H={\cal O}_H(a)\oplus{\cal O}_H(b).$
By \cite[p. 42]{oss} we know that ${\cal E}$ splits if and only if ${\cal E}|H$ splits for some
hyperplane $H$. Hence
$${\cal E}={\cal O}(a)\oplus{\cal O}(b).$$
\end{proof}

\begin{thm}\label{g<16}Let $X\subset{\bb P}({\cal E})$ be a Calabi-Yau-manifold, 
with ${\cal E}\seq{\bb P}^3$ being a rank-2-bundle. For the generic splitting type $(a,b)$ (with $a\le b$) 
of ${\cal E}$ holds $b-a\le 4$ and hence
$$\gamma({\cal E})\le 16,$$
or equivalently,
$$c_3(X)\ge -296.$$
\end{thm}

\begin{proof}Let $X=\{s=0\}$, where $s\in H^0(-K_{{\bb P}({\cal E})})$, $L\subset{\bb P}^3$ a general line
and
$${\cal E}|L={\cal O}(a)\oplus{\cal O}(b),$$
with $a\le b$.
Then $s$ induces a section 
$$t=p_*s\in H^0(S^2{\cal E}\otimes\det{\cal E}^{-1}\otimes{\cal O}(4)).$$
We have
$$S^2{\cal E}\otimes\det{\cal E}^{-1}\otimes{\cal O}(4)|L\cong{\cal O}(a-b+4)\oplus{\cal O}(4)\oplus{\cal O}(b-a+4)$$
and by the general choice of $L$ the intersection $X\cap p^*L$ can be assumed to be smooth. 
Let $d_{00}:=a-b+4, d_{01}:=4, d_{11}:=b-a+4.$
We denote by $[x_0:x_1]$ the coordinates of the fibres of $p$ in a trivialzing neighbourhood
$p^{-1}(U)$, with $U\subset{\bb P}^3$. In this neighbourhood we can express
$s|p^*L$ as
$$s|p^*L=\sum s_{ij}x_ix_j,$$
with $s_{ij}\in H^0({\cal O}(d_{ij}))|U$.
If $b-a>4$, then $H^0({\cal O}(d_{00}))=0$ and hence $s_{00}=0.$ Therefore
$$s|p^*L=x_1(s_{01}x_0+s_{11}x_1)$$
is reducible and therefore $X\cap p^*L$ is singular along
$S:=\{x_1=0\}\cap\{s_{01}x_0+s_{11}x_1=0\}=\{x_1=0\}\cap\{s_{01}=0\}$.
Since $s_{01}\in H^0({\cal O}(4))$, we conclude $S\not=\emptyset$. 
Hence this case does not occur and therefore $b-a\le 4$.

By applying Lemma \ref{3maxg} we conclude finally
$$\gamma({\cal E})\le\gamma({\cal O}\oplus{\cal O}(b-a))\le\gamma({\cal O}\oplus{\cal O}(4))=16.$$
\end{proof}

\subsection{The discriminant map}

To be able to compute $K(X)$ we need information about the morphism $\pi:X\seq{\bb P}^3$, which
is the natural projection.

\begin{constr}\label{3discr}
Let $X=\{s=0\}$, with $s\in H^0(-K_{{\bb P}({\cal E})})$ and
$$V:=\{p\in{\bb P}^3 | \pi \mbox{ is locally in } p\mbox{ not an \'etale covering} \}.$$
We define the discriminant
$$\Delta_{({\cal E},{\cal F})}:S^2{\cal E}\otimes{\cal F}\seq (\det({\cal E})\otimes{\cal F})^{\otimes 2}$$
by
$$\Delta_{({\cal E},{\cal F})}(\sum_{1\le i<j\le 2}c_{ij}s_is_j\otimes f):=(c_{12}^2-4c_{11}c_{22})(s_1\wedge s_2\otimes f)^{\otimes 2},$$
where $s_1,s_2$ is a ${\cal O}(U)$-basis of ${\cal E}(U)$, ${\cal F}\seq{\bb P}^3$ a line bundle and 
$f\in {\cal F}(U)$ a generator of ${\cal F}(U)$ for a small open set $U\subset{\bb P}^3$.

It is an easy computation that this definition is independent of the chosen bases.

Now we specify 
$${\cal F}=\det{\cal E}^\vee\otimes{\cal O}(4).$$
Then the discrimnant is a map
$$\Delta_{\cal E}:p_*(-K_{{\bb P}({\cal E})})\seq {\cal O}(8),$$
with
$$\{\Delta_{\cal E}(p_*s)=0\}=V$$
set theoretically: in local coordinates
$$s=\sum s_{ij}x_ix_j,$$
where $[x_0:x_1]$ denotes the coordinates of the fibre, $V$ is the locus, where the zeroes of
\begin{eqnarray*}&\sum s_{ij}(z)x_ix_j=0&\label{qu}\end{eqnarray*}
are not two distinct points. By definition this is the
discriminant locus of the qudratic equation in $x_0,x_1$, given by
$$s_{01}^2-4s_{00}s_{11}=0.$$
This coincides with the discriminant locus of $p_*s$.

Since on a trivializing neighbourhood $U\subset{\bb P}^3$ the map is given by
$$\Delta_{\cal E}(t)|U=t_{12}^2-4t_{11}t_{22},$$
if $t\in H^0(p_*(-K_{{\bb P}({\cal E})}))$ and $t|U=(t_{11},t_{12},t_{22})$,
we see, that, in particular, $H^0(\Delta_{\cal E})$ is a holomorphic map.


Moreover, 
$$H^0(\Delta_{\cal E})(rt)=r^2H^0(\Delta_{\cal E})(t)$$
for $r\in{\bb C}, t\in H^0(-K_{{\bb P}({\cal E})})$. 
Hence we can projectivize, but cannot exclude that
$H^0(\Delta_{\cal E})(s')=0$ for some $s'\not= 0$. Therefore we get a rational map
$$\delta_{\cal E}:{\bb P}(H^0(-K_{{\bb P}({\cal E})}))\ratmap{\bb P}(H^0({\cal O}(8)))\cong{\bb P}^{164}.$$

Let for the moment $V':=\{z\in{\bb P}^3| H^0(\Delta_{\cal E})(s)(z)=0\}$ in the sense of ideals.
If we denote
$$P:=\{z\in{\bb P}^3 | \dim \pi^{-1}(z)=1\},$$
then we see, that
$$P=\{z\in{\bb P}^3| \sum s_{ij}x_ix_j=0 \mbox{ f\"ur alle}[x_o:x_1]\}$$
and hence
$$P=\{z\in{\bb P}^3| s_{00}(z)=s_{01}(z)=s_{11}(z)=0\}\subset Sing(V').$$
Moreover, this shows
$$P=\{z\in{\bb P}^3 | \pi^{-1}(z)\cong{\bb P}^1\}.$$

Now let $z\in Sing(V')$. If $s_{00}(z)=s_{01}(z)=s_{11}(z)=0$, then $z\in P$. 
So let us assume $s_{00}(z)\not= 0$ or $s_{01}(z)\not= 0$. Let us define
$$x:=[s_{01}(z):-2s_{00}(z)]\in p^{-1}(z).$$
Since $z\in V$, we get that $\Delta_{\cal E}(s)(z)=s_{01}(z)^2-4s_{00}(z)s_{11}(z)=0.$ Therefore
$$s(x)=s_{00}(z)s_{01}(z)^2-2s_{00}(z)s_{01}(z)^2+4s_{11}(z)s_{00}(z)^2=-s_{00}(z)\Delta_{\cal E}(s)(z)=0,$$
hence $x\in X$. 

We want to show that $x\in X$ is singular. For this we have to compute in the point $x$
\begin{eqnarray}
\frac{\partial s}{\partial x_0}=&2s_{00}x_0+s_{01}x_1&=0\label{3x1}\\
\frac{\partial s}{\partial x_1}=&s_{01}x_0+2s_{11}x_1&=0\label{3x2}\\
\frac{\partial s}{\partial z_i}=&\frac{\partial s_{00}}{\partial z_i}x_0^2+\frac{\partial s_{01}}{\partial z_i}x_0x_1+
\frac{\partial s_{11}}{\partial z_i}x_1^2&=0\label{3x3}
\end{eqnarray}
and we know moreover, since $z\in Sing(V')$, that in the point $z$ holds
\begin{eqnarray}
s_{01}^2-4s_{00}s_{11}&=&0\label{3z1}\\
2s_{01}\frac{\partial s_{01}}{\partial z_i}-4s_{11}\frac{\partial s_{00}}{\partial z_i}-4s_{00}\frac{\partial s_{11}}{\partial z_i}&=&0.
\label{3z2}
\end{eqnarray}

Using the expression for  $x$ in (\ref{3x1}), (\ref{3x2}) and (\ref{3x3}) we compute
\begin{eqnarray}
\frac{\partial s}{\partial x_0}=&2s_{00}s_{01}-2s_{00}s_{01}&=0\nonumber\\
\frac{\partial s}{\partial x_1}=&s_{01}^2-4s_{00}s_{11}&=0\nonumber\\
\frac{\partial s}{\partial z_i}=&\frac{\partial s_{00}}{\partial z_i}s_{01}^2-2\frac{\partial s_{01}}{\partial z_i}s_{00}s_{01}+
4\frac{\partial s_{11}}{\partial z_i}s_{00}^2&=\nonumber\\
=&4\frac{\partial s_{00}}{\partial z_i}s_{00}s_{11}-2\frac{\partial s_{01}}{\partial z_i}s_{00}s_{01}+
4\frac{\partial s_{11}}{\partial z_i}s_{00}^2&=\nonumber\\
=&-s_{00}(2s_{01}\frac{\partial s_{01}}{\partial z_i}-4s_{11}\frac{\partial s_{00}}{\partial z_i}-4s_{00}\frac{\partial s_{11}}{\partial z_i})&=0,\nonumber
\end{eqnarray}
with the lase equation using (\ref{3z1}) as well as (\ref{3z2}).

Thus we have proved that $x\in X$ is singular. But we assumed $X$ to be smooth. Hence it is proven that
$P=Sing(V').$
In particular, $V'$ is reduced and therefore $V'=V$ in the sense of ideals.

Now we know $V=\delta_{\cal E}(X)\in|{\cal O}(8)|$ and 
$$P=Sing(V).$$
  
If we assume $\rho(X)=2$, we will see in the next section that in this case $P\not=\emptyset$.
The image of $\delta_{\cal E}$ then is a subvariety of the singular octics in ${\bb P}^3$, in
particular 
$\delta_{\cal E}$ is not surjective.
\end{constr}

At this point we should mention the work of Clemens, Cynk and Szemberg \cite{cl83,cs99,cy99}. The first
article mentioned describes the construction of Calabi-Yau-threefolds as resolutions of
double covers of ${\bb P}^3$ ramified over a given (singular) octic. 
The latter two papers are concerned with the euler number of such Calabi-Yau-threefolds. The track followed
in the present paper reverses the direction, since we are given first the Calabi-Yau-threefold and then
construct the octic. Hence our method can be used to construct octic hypersurfaces with
many nodes. This will be done in another paper.

To prove now the finiteness of $P$ we will fall back upon the bounds for 
$\gamma$ proved in the last section.

\begin{lemma}\label{dimP=0}Let $X\subset{\bb P}({\cal E})=:Z$ be a Calabi-Yau manifold, 
with ${\cal E}\seq{\bb P}^3$ being a rank-2-bundle. Furthermore let
$P:=\{p\in{\bb P}^3 | \pi^{-1}(p)\cong{\bb P}^1\}.$ Then $\dim P\le 0$ holds true.
\end{lemma}

\begin{proof}Let us assume first $\dim P=2$. Then
$\dim \pi^* P=3$ ans since $\pi^* P\not= X$ we conclude, that $X$ is reducible, hence not smooth.

If $\dim P=1$, then $D:=\pi^* P\in Pic(X)$ is an effective divisor, satisfying
$$D.\pi^*h^2=0.$$ 
If we write
$$\mu D={\cal O}_X(2)+k\pi^*h,$$
$\mu\in{\bb R}^+$, we compute
\begin{eqnarray*}0&=&D.\pi^*h^2\\
&=&({\cal O}_Z(2)+kp^*h).p^*h^2.({\cal O}_Z(2)+(4-c_1({\cal E}).h^2)p^*h\\
&=&8+2k+2c_1({\cal E}).h^2.\end{eqnarray*}
Hence 
$$\mu D=-K_Z-8p^*h|X.$$
But we know additionally
\begin{eqnarray*}0&<&\deg P\\
&=&D.\pi^*h.{\cal O}_X(1)\\
&=&\frac{1}{\mu}(-K_Z-8p^*h).\pi^*h.{\cal O}_Z(1).(-K_Z)\\
&=&\frac{1}{\mu}(\gamma-16)\le 0,\end{eqnarray*}
what is a contradiction.
\end{proof}

This we now use to prove a result about the number of full fibres in $X$:

\begin{thm}\label{3g<16}Let $X\subset{\bb P}({\cal E})=:Z$ be a Calabi-Yau-manifold, 
with ${\cal E}\seq{\bb P}^3$ being a rank-2-bundle. 
$X$ contains exactly $64-4\gamma$ fibres of $p$.
\end{thm}

\begin{proof}
Let $s\in H^0(-K_Z)$ be such that $X=\{s=0\}$. 
This section $s$ induces a section $t\in H^0(p_*(-K_Z))$. Let
$$P:=\{t=0\}.$$
Since by Lemma \ref{dimP=0} holds $\dim P\le 0$, we know
$$[P]=c_3(p_*(-K_Z)).$$.

To compute $c_3(p_*(-K_Z))$ we make the usual ansatz
$$c({\cal E})=(1+at)(1+bt).$$
Then we can express
$$c(S^2{\cal E}\otimes{\cal O}(r))=(1+(2a+r)t)(1+(a+b+r)t)(1+(2b+r)t).$$
If we do the multiplications and replace again by Chern classes of ${\cal E}$, we get
$$c_3(S^2{\cal E}\otimes{\cal O}(r))=4c_2({\cal E})c_1({\cal E})+2r(c_1^2({\cal E}).h+2c_2({\cal E}).h)+3r^2c_1({\cal E})
+r^3.$$
Setting $r=4-c_1({\cal E}).h^2$ finally leads to
$$c_3(p_*(-K_Z))=64-4\gamma.$$
\end{proof}

\subsection{The generalized Koll\'ar Lemma}

By Theorem \ref{3g<16} the only possible case, in which $\pi^*h$ is ample, can be $\gamma=16$. But
Lemma \ref{3maxg} and Lemma \ref{g<16} state, that in this case, ${\cal E}={\cal O}\oplus{\cal O}(4)$.
We computed in Example \ref{3rho=1} that then $\rho(X)=1$. Hence we proved

\begin{cor}\label{hnef}Let $X\subset{\bb P}({\cal E})=:Z$ be a Calabi-Yau-manifold, 
with ${\cal E}\seq{\bb P}^3$ being a rank-2-bundle. Assume that $\rho(X)=2$. Then $\pi^*h\in\partial K(X)$.
\end{cor}

So we turn our attention to the 'other' side of the ample cone. The following arguments are similar to
the corresponding section of \cite{k01}.

\begin{lemma}\label{3lemkol}Let $X\subset{\bb P}({\cal E})=:Z$ be a Calabi-Yau-manifold with $\rho(X)=2$, 
with ${\cal E}\seq{\bb P}^3$ being a rank-2-bundle. If $-K_Z$ is not nef, then
$K(X)=K(Z)|X$.
\end{lemma}

\begin{proof}Let $D\in Pic(Z)\otimes{\bb Q}$ such that  $D|X$ is nef. 
Without loss of generality we may assume $D={\cal O}_Z(2)+kp^*h, k\in{\bb Q}$.
Let us denote $E:=D+K_Z$. Then $E=lp^*h,l\in{\bb Q}$, hence $E$ nef oder $-E$ nef. If $-E$ is nef, then
$-K_Z$ is also nef, hence we get a contradiction. Therefore $E$ is nef. Now let $C\subset Z$ be a curve. 
If $C\subset X$, then by assumption $D.C\ge 0$. But if $C\not\subset X$, then $-K_Z.C\ge 0$ since
$X\in|-K_Z|$ and hence
$$D.C=(E-K_Z).C>0.$$
This shows that $D$ is nef.
\end{proof}

The generalization of the Lemma of Koll\'ar is possible in an unrestricted way:

\begin{thm}\label{3genkol2}Let $X\subset{\bb P}({\cal E})=:Z$ be a Calabi-Yau-manifold with $\rho(X)=2$, 
with ${\cal E}\seq{\bb P}^3$ being a rank-2-bundle.
Then holds
$$K(X)=K(Z)|X.$$
\end{thm}

We divide the proof in several steps:

\begin{lemma}\label{smrat}
Let $Z$ be a fourfold, such that $-K_Z$ is big and nef, but not ample. 
If $\Phi_{|-mK_Z|}:Z\seq Z'$ contracts only a finite number of curves,
then these are smooth and rational.
\end{lemma}

\begin{proof}
Let $C$ be an irreducible curve, contracted by $\Phi$.
Since $-mK_Z=\Phi^*{\cal O}_{Z'}(1)$, also $-mK_{Z'}={\cal O}_{Z'}(1)$ holds and $Z'$ 
has only canonical singularities, in particular, they are rational.
This tells us $R^1\Phi_*{\cal O}_{Z}=0$ and by
$$0\seq{\cal I}_C\seq{\cal O}_Z\seq{\cal O}_C\seq 0$$
we conclude $R^1\Phi_*{\cal O}_C=0$. By the Leray spectral sequence follows $H^1({\cal O}_C)=0$. 
Hence $C$ is smooth and rational. 
\end{proof} 

\begin{lemma}\label{3nocurves}Let $Z$ be a fourfold such that $-K_Z$ is  big and nef, but not ample. 
Then the exceptional locus of $\Phi_{|-mK_Z|}:Z\seq Z'$ contains a two-dimensional component.\end{lemma}

\begin{proof}
Let us assume, $\Phi$ contracts only a finite number of curves. By Lemma \ref{smrat} these are
smooth and rational. Let $C$ be such a curve.
ABy the adjunction formula $K_Z.C=0$ implies
$c_1(N_{C|Z})=c_1(K_C)=-2$. Now we compute
$$\chi(N_{C|Z})=3(1-g(C))+c_1(N_{C|Z})=1>0.$$
Therefore $C$ deforms in $Z$, hence there is a surface contracted, contradicting the assumption.
\end{proof}

\begin{lemma}\label{3generators}
Let $X\subset{\bb P}({\cal E})=:Z$ be a Calabi-Yau-manifold with $\rho(X)=2$, 
with ${\cal E}\seq{\bb P}^3$ being a rank-2-bundle.
Then 
$$H^4(Z,{\bb Z})=<{\cal O}_Z(1).p^*h,(p^*h)^2>.$$\end{lemma}

\begin{proof}By the K\"unneth formula $b_4(Z)=2$.

Tos show that $v_1:=(p^*h)^2$ und $v_2:={\cal O}_Z(1).p^*h$ form a ${\bb Z}$-Basis for $H^4(Z,{\bb Z})$, 
it suffices to show that the matrix $A=(a_{ij})=(v_i.v_j)$ is invertible over ${\bb Z}$. 
But by Lemma \ref{3intersection}
$$A=\left( \begin{array}{ccc}
0 & 1\\
1 & c_1({\cal E}).h^2
\end{array} \right)$$
what proves the Lemma.
\end{proof}

\begin{lemma}\label{3noniso}
Let $X\subset{\bb P}({\cal E})=:Z$ be a Calabi-Yau-manifold with $\rho(X)=2$, 
with ${\cal E}\seq{\bb P}^3$ being a rank-2-bundle.
Assume that $-K_Z$ is big and nef, but not ample. Then
$\Phi_{|-mK_Z|}|X:X\seq X'$ is no isomorphism.
\end{lemma}

\begin{proof}Let $E\seq V$ be the exceptional locus of $\Phi$. 
Since 
$$kX=\phi^*H,$$
for a $k\in{\bb Z}$.
Since $H$ is ample, $H$ intersects with every positive dimensional component of $V$. This implies
that $\phi|X$ can be an isomorphism only if $\dim V=0$.

So let us assume $\dim V=0$.
By Lemma \ref{3nocurves} there is a surface $G\subset Z$, which gets contracted to a point
by $\Phi$. In particular, $-K_Z.G\equiv 0$. Let $a,b\in{\bb Z}$ such, that
$$[G]=a{\cal O}_Z(1).p^*h+b(p^*h)^2.$$
Then $-K_Z.G.p^*h=0$ and $-K_Z.G.{\cal O}_Z(1)=0$ are equivalent to
\begin{eqnarray}
(c_1({\cal E}).h^2+4)a+2b&=&0\nonumber\\
(c_1^2({\cal E}).h-2c_2({\cal E}).h+4c_1({\cal E}).h^2)a+(c_1({\cal E}).h^2+4)b&=&0,\nonumber
\end{eqnarray}
what has a non-trivial solution only if
$$\det \left( \begin{array}{cc} (c_1({\cal E}).h^2+4)& 2\\
 (c_1^2({\cal E}).h-2c_2({\cal E}).h+4c_1({\cal E}).h^2) &(c_1({\cal E}).h^2+4)\\
\end{array} \right)  =0.$$

This means exactly $\gamma=16$. 
As argued above, this amounts to $\rho(X)=1$, hence contradicts to our assumption.
\end{proof}

Now we can give the proof of Theorem \ref{3genkol2}:

\begin{proof}[of Theorem \ref{3genkol2}]
We assume $K(X)\not=K(Z)|X$. The Koll\'ar lemma says that then $-K_Z$ is not ample. But
by Corollary \ref{hnef} the pullbacks $p^*h$ and $\pi^*h$ nef and not ample. 
So, if $K(X)\not= K(Z)|X$,
by Lemma \ref{3lemkol} the divisor $-K_Z|X$ has to be ample und $-K_Z$ has to be big and nef, but not ample.
By 
$$0\seq -(m-1)K_Z\seq -mK_Z\seq -mK_Z|X\seq 0$$
and $H^1(-(m-1)K_Z)=0$ for $m>0$ we see that
$$H^0(-mK_Z)\seq H^0(-mK_Z|X)$$ is surjective, hence $\phi_{|-mK_Z|_X|}=\phi_{|-mK_Z|}|X$. 
Since $-K_Z|X$ is ample, this means that $\phi_{|-mK_Z|}|X$ is an isomorphism for $m\gg 0$. 
This contradicts to
Lemma \ref{3noniso}.
\end{proof}

\subsection{Rationality of $K(X)$ and Positivity of $c_2(X)$}

\begin{cor}\label{3c2g0}Let $X\subset{\bb P}({\cal E})=:Z$ be a Calabi-Yau-manifold with $\rho(X)=2$, 
with ${\cal E}\seq{\bb P}^3$ being a rank-2-bundle.
Then
$$D.c_2(X)>0 \mbox{ for all } D\in\overline{K(X)}.$$
\end{cor}

\begin{proof}If $-K_Z$ is not ample, then by Theorem \ref{3genkol2} also $-K_Z|X$ is not ample. 
Hence there is some non-negative $k\in{\bb Q}$, such that
$$D:=-K_Z|X+k\pi^*h \in \partial K(X).$$
By Lemma \ref{3c2} holds 
$$D.c_2(X)=-K_Z|X.c_2(X)+k\pi^*h.c_2(X)\ge 56+44k>0.$$
By Theorem \ref{3genkol2} follows, that $\pi^*h\in\partial K(X)$ and by  $\pi^*h.c_2(X)=44$ 
the claim follows.

If $-K_Z$ is ample, the claim of the theorem has been proven by Oguiso and Peternell in \cite{op}.
\end{proof}
 
The rationality of the K\"ahler cone follows also easily under mild restrictions:

\begin{thm}
Let $X\subset{\bb P}({\cal E})=:Z$ be a Calabi-Yau-manifold with $\rho(X)=2$, 
with ${\cal E}\seq{\bb P}^3$ being a rank-2-bundle satisfying
$h^0(-K_{{\bb P}({\cal E})})>1$. Then $\partial K(X)$ is rational. 
\end{thm}

\begin{proof}With respect to \ref{3genkol2} it suffices to show, that 
$\partial K({\bb P}({\cal E}))$ is rational.
Denote like before $Z:={\bb P}({\cal E})$. 

We consider three cases: $-K_Z$ is ample, $-K_Z$ is nef, but not ample and finally, $-K_Z$ is not nef.

Let $-K_Z$ be ample. Then the cone theorem states the rationality of $\partial K(Z)$. 

Now let $-K_Z$ be nef, but not ample. Then $-K_Z\in\partial K(X)$ and hence $\partial K(X)$ is rational.

Finally, let $-K_Z$ be not nef. By Theorem \ref{3genkol2} holds 
$$K(X)=K(Z)|X.$$
Since $h^0(-K_Z)>1$ holds, $-K_Z|X$ is effective. By assumption also  $-K_Z|X$ is not nef. 
Now we want to use
the log-rationality theorem (\cite[Thm 4-1-1]{kmm}), which states that
$$\sup\{r\in\bb R|H+r(K_X+\Delta)\in K(X)\}\in\bb Q,$$
if $H\in Pic(X)$ is ample, $K_X+\Delta$ not nef and $\Delta$ an effective ${\bb Q}$-Divisor, 
such that $(X,\Delta)$ has only weak log-terminal singularities. The latter property can be
reached by choosing $\varepsilon\Delta$ instead of $\Delta$ for $0<\varepsilon\ll 1$.
Note that $K_X=0$.

Because $-K_Z|X$ is not nef, we apply
the log-rationality theorem for an arbitrary ample $H\in Pic(X)$ and 
$\Delta:=\varepsilon(-K_Z|X)$ for 
$0<\varepsilon\ll 1$. In this way we get that $\partial K(X)$ is rational.
\end{proof}

A result of
Wilson \cite[Facts A,B,C]{wi94} says, that the rationality of $\partial K(X)$ und the positivity
of $c_2(X)$ imply, that 
$\partial K(X)$ is semiample. Hence we get

\begin{cor}
Let $X\subset{\bb P}({\cal E})=:Z$ be a Calabi-Yau-manifold with $\rho(X)=2$, 
with ${\cal E}\seq{\bb P}^3$ being a rank-2-bundle. If $\partial K(X)$ is rational,
then there is a second contraction $X\seq X'$ (apart from the first Stein factor of $\pi$).
\end{cor}

\section{Calabi-Yau-threefolds in ${\bb P}^3$-bundles over ${\bb P}^1$}

According to Corollary \ref{p1}, if $X$ is a Calabi-Yau-threefold in a ${\bb P}^3$-bundles over ${\bb P}^1$
with $\rho(X)=2$, the K\"ahler cone is rational. 
Moreover, with similar arguments as in the previous section, it can be 
proved, that $K(X)=K(Z)|X$. 

We meet the convention, that ${\cal E}\seq{\bb P}^1$ is normalized in such a way, that
$${\cal E}={\cal O}\oplus{\cal O}(a_1)\oplus{\cal O}(a_2)\oplus{\cal O}(a_3)$$
with $a_1\le a_2\le a_3$ and $X\subset{\bb P}({\cal E})=:Z$ is the considered Calabi-Yau-manifold.

\begin{thm}\label{4genkol}
Let $X\subset{\bb P}({\cal E})=:Z$ be a Calabi-Yau-threefold with ${\cal E}\seq{\bb P}^1$
a rank-4-bundle. Then $K(X)=K(Z)|X$. 
\end{thm}

\begin{proof}We consider ${\cal E}$ normalized as above. In particular, ${\cal O}_Z(1)\in\partial K(Z)$.
The canonical bundle of $Z$ is
$$-K_Z={\cal O}_Z(4)+(2-c_1({\cal E}))p^*h$$
and $X\in|-K_Z|$. 

If $c_1({\cal E})>2$, then
$$-K_Z.{\bb P}({\cal O})<0,$$
hence ${\bb P}({\cal O})\subset X$. Therefore ${\cal O}_X(1)$ is also nef and not ample.

If $c_1({\cal E})<2$, then $-K_Z$ is ample and we use the Lemma of Koll\'ar (Lemma \ref{0kollar}).

If $c_1({\cal E})=2$, then $a_1=0$ follows and therefore ${\bb P}({\cal O})$ deforms in $Z$. Since
$-K_Z.{\bb P}({\cal O})=0$,
either there is a curve of the form ${\bb P}({\cal O})$ lying in $X$ or the surface 
$G:={\bb P}({\cal O}\oplus{\cal O})$ has the property
$$X\cap G=\emptyset.$$
In the first case $K(X)=K(Z)|X$ is proven.

In the second case we compute like in the previous section
\begin{eqnarray}&<{\cal O}_Z(1).p^*h, {\cal O}_Z(1)^2>_{\bb Z}=H^4(Z,{\bb Z}).&\label{4gen}\end{eqnarray}
Again, we show this by setting $v_1:={\cal O}_Z(1).p^*h, v_2:={\cal O}_Z(1)^2$ and computing
the matrix $A=(v_iv_j)_{ij}$:
$$A=\left(\begin{array}{cc}0 & 1\\ 
1 & c_1({\cal E}) \end{array}\right).$$
Obviously, $A\in Gl(2,{\bb Z})$ and therefore claim (\ref{4gen}) is proven.

Now we can write
$$[G]=kv_1+lv_2.$$
With this we see
$$0=-K_Z.G.p^*h=4l$$
and hence
$$0=-K_Z.G.{\cal O}_Z(1)=4k+4lc_1({\cal E})+(2-c_1({\cal E}))l=4k,$$
therefore
$$[G]=0,$$
what is a contradiction to the projectivity of $Z$.
\end{proof}

The intersection theory on $X$ looks like follows:

\begin{enumerate}
\item $c_3(X)=-168$
\item ${\cal O}_Z(1)|X.c_2(X)=6c_1({\cal E})+44, \pi^*h.c_2(X)=24$, hence $-K_Z|X.c_2(X)=224.$ 
\item $(-K_Z|X)^3=512, (-K_Z|X)^2.\pi^*h=64$
\end{enumerate}

With this, we already are able to prove the positivity of $c_2(X)$:

\begin{thm}Let $X\subset{\bb P}({\cal E})=:Z$ be a Calabi-Yau-threefold with ${\cal E}\seq{\bb P}^1$
a rank-4-bundle normalized as above. Then
$$D.c_2(X)>0$$
for all $D\in\partial K(X)$.
\end{thm}

\begin{proof}By the chosen normalization, ${\cal O}_Z(1)$ is nef and not ample,
hence Theorem \ref{4genkol} shows that ${\cal O}_X(1)$ is also nef and not ample. 
Since $c_1({\cal E})\ge 0$,
the above formulas imply ${\cal O}_X(1).c_2(X)>0$ and $\pi^*h.c_2(X)>0$. Hence
the claim is proved.
\end{proof}

Again we refer to
 
\begin{eqnarray}
&0\seq T_{Z|{\bb P}^1}\seq T_Z \seq p^* T_{{\bb P}^1}\seq 0&\label{f1}\\
&0\seq {\cal O}_Z\seq p^*({\cal E}^\vee)\otimes{\cal O}_Z(1) \seq T_{Z|{\bb P}^1}\seq 0&\label{f2}\\
&0\seq T_X\seq T_Z|X \seq N_{X|Z}\seq 0&\label{f3}
\end{eqnarray}

For the later description we compute the Picard number.

\begin{thm}Let $X\subset{\bb P}({\cal E})=:Z$ be a Calabi-Yau-threefold with ${\cal E}\seq{\bb P}^1$
a rank-4-bundle normalized as above. Then
$$\rho(X)=2+h^1(-K_Z).$$
In particular,
$$\rho(X)=2\iff c_1({\cal E}\le 3.$$
\end{thm}

\begin{proof}
We look at the two sequences
\begin{eqnarray}&0\seq N_{X|Z}^\vee\seq\Omega_Z|X\seq\Omega_X\seq 0&\label{?}\end{eqnarray}
and
\begin{eqnarray}&0\seq\Omega_Z\otimes K_Z\seq\Omega_Z\seq\Omega_Z|X\seq 0.&\label{??}\end{eqnarray}

Our first aim is to show $H^i(T_Z)=H^i({\cal E}^\vee\otimes{\cal E})$ for $i>1$. For this purpose we calculate
$R^ip_*(p^*({\cal E}^\vee)\otimes{\cal O}_Z(1))={\cal E}^\vee\otimes R^ip_*{\cal O}_Z(1)=0$ for $i>0$. Therefore the
Leray spectral sequence yields 
$$H^i(p^*({\cal E}^\vee)\otimes{\cal O}_Z(1))=H^i({\cal E}^\vee\otimes{\cal E}).$$
Sequence (\ref{f2}) shows, 
that $H^i(p^*({\cal E}^\vee)\otimes{\cal O}_Z(1))=H^i(T_{Z|{\bb P}^2})$ for $i>0$,
since by $R^ip_*{\cal O}_Z=0$ for $i>0$ we get $H^i({\cal O}_Z)=H^i({\cal O}_{{\bb P}^2})=0$ for $i>0$.
For applying (\ref{f1}), we verify by the projection formula
$R^ip_*p^*T_{{\bb P}^2}=0$ for $i>0$. 
Hence again the Leray spectral sequence implies
$$H^i(p^*T_{{\bb P}^2})=H^i(T_{{\bb P}^2})=0$$
for $i>0$. This implies with (\ref{f1}) now, that
$$H^i(T_{Z|{\bb P}^2})=H^i(T_Z)$$
for $i>1$. Therefore
$$H^i(T_Z)=H^i({\cal E}^\vee\otimes{\cal E})$$
for $i>1$.

So we see that
$$H^i(\Omega_Z\otimes K_Z)=H^{4-i}(T_Z)^\vee=H^{4-i}({\cal E}^\vee\otimes{\cal E})^\vee$$
for $i<3$. In particular,
$$H^1(\Omega_Z\otimes K_Z)=H^3({\cal E}^\vee\otimes{\cal E})^\vee=0.$$
Since $H^{2,1}(Z)=0, H^2({\cal E}^\vee\otimes{\cal E})=0$ and 
$N_{X|Z}^\vee=K_Z|X$ the cohomology sequences of (\ref{?}) and (\ref{??}) contain 
\begin{eqnarray}&0\seq H^1(\Omega_Z)\seq H^1(\Omega_Z|X)\seq 0&\label{!4}
\end{eqnarray}
resp.\,
\begin{eqnarray}&0\seq H^1(K_Z|X)\seq H^1(\Omega_Z|X)\seq H^1(\Omega_X)\seq H^2(K_Z|X).&\label{!!4}\end{eqnarray}

The sequence
$$0\seq {\cal O}_Z\seq -K_Z\seq -K_Z|X\seq 0$$
and Serre duality imply
$$h^i(-K_Z)=h^i(-K_Z|X)=h^{3-i}(K_Z|X)$$
for $i>0$. Since $R^ip_*(-K_Z)=0$ for $i>0$ we compute
$$h^2(-K_Z)=h^2(p_*(-K_Z))=0.$$

Considering the sequences (\ref{!4}) and (\ref{!!4}) together with this vanishing yields
$$\rho(X)=2+h^2(K_Z|X)=2+h^1(-K_Z).$$
Since $R^1p_*(-K_Z)=0$, we compute further
$$h^1(-K_Z)=h^1(S^4{\cal E}\otimes{\cal O}(2-c_1({\cal E}))).$$
Hence the condition $\rho(X)=2$ is equivalent to
$$2-c_1({\cal E})\ge -1,$$
what proves the claim. 
\end{proof}

In the chosen normalization of ${\cal E}$ the divisor ${\cal O}_Z(1)$ is nef, but not ample. By 
$K(X)=K(Z)|X$ also ${\cal O}_X(1)$ is nef and not ample. The computation
$${\cal O}_X(1)^3={\cal O}_Z(1)^3.(-K_Z)=3c_1({\cal E})+2>0$$
shows with the base point free theorem, that there is a birational contraction
$$\psi:X\seq X'.$$

This contraction is described as follows:

\begin{thm}
\label{exz}Let $X\subset{\bb P}({\cal E}):=Z$ be a Calabi-Yau-manifold with $\rho(X)=2$, with
${\cal E}\seq{\bb P}^1$ being a rank-4-bundle normalized as above.
Let $\psi:X\seq X'$ be the birational contraction and $E$ its exceptional locus.
Denote by ${\cal F}:={\cal O}\bigoplus_{i|a_i=0}{\cal O}(a_i)$ the maximal trivial subbundle of ${\cal E}$. 
\begin{enumerate}
\item If $c_1({\cal E})=3$, then $rk{\cal F}\le 2$ and $E={\bb P}({\cal F})\cong{\bb P}^1\times{\bb P}^{rk{\cal F}-1}$.
\item If $c_1({\cal E})=2$, then $rk{\cal F}\in\{2,3\}$ and $E=X\cap{\bb P}({\cal F})={\bb P}^1\times Y$, 
with $\dim Y=rk{\cal F}-2.$ 
\begin{enumerate}
\item If $rk{\cal F}=2$, then $Y$ consists of four points.
\item If $rk{\cal F}=3$, then $Y$ is a smooth plane quartic curve.
\end{enumerate}
\item If $c_1({\cal E})=1$, then $Z$ is the blowup of ${\bb P}^4$ in a linearly embedded ${\bb P}^2$; if 
$X\in|-K_Z|$ is chosen generally, then $E=\bigcup_{i=1}^{16} C_i$, with $C_i\cong{\bb P}^1$; in this
case, $X'$ is a quintic in ${\bb P}^4$ with $16$ double points on a lineraly embedded ${\bb P}^2$.
\item If $c_1({\cal E})=0$, then $Z={\bb P}^1\times{\bb P}^3$ and $E=\bigcup_{i=1}^{64}C_i,$ with
$C_i\cong{\bb P}^1$.
\end{enumerate}
\end{thm}

\begin{proof}${\cal E}$ is normalized in such a way that ${\cal O}_Z(1)\in\partial K(Z)$. 
First, we show that all curves $C$, which satisfy ${\cal O}_Z(1).C=0$, have the form
${\bb P}({\cal O})$ for some quotient ${\cal E}\seq{\cal O}\seq 0$:
If ${\cal O}_Z(1).C=0$, then $s|C$ is constant for all $s\in H^0({\cal O}_Z(1))$. Writing
$$s=\sum s_ix_i,$$
with $[x_0:...:x_3]$ being fibre coordinates, $[z_0:z_1]$ coordinates of the  base-${\bb P}^1$, 
and $s_i\in H^0({\cal O}(a_i))$,
constancy of $s|C$ for all $s\in H^0({\cal O}_Z(1))$ implies, that
$$x_i=0 \mbox{, if }a_i>0 \mbox{ and } x_i=c_i \mbox{, if }a_i=0$$
for $c_i\in{\bb C}$. This proves the claim.

$(i)$: In case $c_1({\cal E})=3$ we compute
for $C={\bb P}({\cal O})$
$$-K_Z.C=({\cal O}_Z(4)+(2-c_1({\cal E}))p^*h).C=2-c_1({\cal E})<0,$$
hence $C\subset X$. Since every curve $C$ satisfying ${\cal O}_Z(1).C=0$ is of the form
${\bb P}({\cal O})$, this implies
$$E={\bb P}({\cal F}).$$
The divisor $E$ is also the exceptional locus of $\tilde\psi:=\Phi_{|{\cal O}_Z(m)|}:Z\seq Z'$ and hence
by Lemma \ref{3nocurves} we know $\dim E\le 2$. Therefore also $rk{\cal F}\le 2$.

$(ii)$: If $c_1({\cal E})=2$, then for $C={\bb P}({\cal O})$
$$-K_Z.C=0,$$
hence $C\subset X$ or $C\cap X=\emptyset$. Again follows
$$E=X\cap{\bb P}({\cal F})={\bb P}^1\times Y.$$
Since $-K_Z={\cal O}_Z(4)$, 
there is an ample hypersurface $V\subset Z'$, such that
$$kX=\tilde\psi^* V,$$
with $\tilde\psi:Z\seq Z'$ like above the birational map given by the linear sytem $|{\cal O}_Z(m)|$ 
In particular,
$$\dim E=\dim{\bb P}({\cal F})-1=rk{\cal F}-1,$$
and hence
$$\dim Y=rk{\cal F}-2.$$

If $rk{\cal F}=2$, the number of contracted curves is computed as
$$Y.F={\cal O}_Z(4).{\bb P}({\cal F}).F=4,$$
since every contracted is of the form ${\bb P}({\cal O})$.

If $rk{\cal F}=3$, i.e. ${\cal E}=3{\cal O}\oplus{\cal O}(2)$, then by \cite{gr2} $Y$ is a smooth
curve, since $\rho(X)=2$ implies that
every contraction is primitive. For a general fibre $F$ if $p$ the intersection $E\cap F$ is
smooth (\cite[Corollary 10.9.1]{ha}) and 
$$\psi|_{E\cap F}:E\cap F\seq Y$$
is an isomorphism, since $F.{\bb P}({\cal O})=1$. 
Furthermore $E\cap F\subset{\bb P}({\cal F})\cap F\cong{\bb P}^2$ and we compute
$$\deg E\cap F=E.F.{\cal O}_Z(1)=-K_Z.{\bb P}({\cal F}).F.{\cal O}_Z(1)=4,$$
because $-K_Z={\cal O}_Z(4)$.

$(iii)$: If $c_1({\cal E})=1$, only ${\cal E}=3{\cal O}\oplus{\cal O}(1)$ is possible and therefore
$Z$ is the blowup of ${\bb P}^4$ in a linearly embedded ${\bb P}^2$. 
Let denote the blowup by $\phi:Z\seq{\bb P}^4$.
Like above the exceptional curves have the form $C={\bb P}({\cal O})$. 
Let $E'$ be the exceptional divisor of the blowup. Then
$$-K_Z=\phi^*{\cal O}_{{\bb P}^4}(5)-E'={\cal O}_Z(5)-E',$$
on the other hand
$$-K_Z={\cal O}_Z(4)+p^*h,$$
hence
$$E'={\cal O}_Z(1)-p^*h.$$
Denote $Y:=X\cap E'$. The sequence
$$0\seq -K_Z-E'\seq -K_Z\seq -K_Z|E'\seq 0$$
implies by $H^1(-K_Z-E')=H^1({\cal O}_Z(3)+2p^*h)=0$ (with Kodaira vanishing), that
$$H^0(-K_Z)\seq H^0(-K_Z|E')$$
is surjective. Since $-K_Z|E'$ is globally generated, 
we see, that by a general choice of $X$ also $Y$ can be assumed to be smooth.

Since ${\cal O}_Z(1).{\bb P}({\cal O})=0$ and $p^*h.{\bb P}({\cal O})=1$ we conclude
$-K_Z.{\bb P}({\cal O})=1$ and therefore
$$\phi|Y:Y\seq{\bb P}^2$$
is generically $1:1$. If $C={\bb P}({\cal O})\subset Y$, then
$$0\seq T_{{\bb P}({\cal O})}\seq T_{Y}|{\bb P}({\cal O})\seq N_{{\bb P}({\cal O})|Y}\seq 0,$$
$T_{{\bb P}({\cal O})}={\cal O}(2)$ as well as 
$$c_1(T_{Y}|{\bb P}({\cal O}))=c_1(N_{Y|X}^\vee|{\bb P}({\cal O}))=-E'.{\bb P}({\cal O})=(-{\cal O}_Z(1)+p^*h).{\bb P}({\cal O})=1,$$
imply that
$$N_{{\bb P}({\cal O})|Y}={\cal O}(-1).$$
This shows that $\phi|Y$ is the blowup of ${\bb P}^2$ in $k$ points. To compute $k$, 
we apply the canonical bundle formula,
$$K_Y^2=(\phi|Y^*({\cal O}_{{\bb P}^2}(3))-\sum_{i=1}^k E_i)^2=9-k,$$
with $E_i$ denoting the exceptional curves. We know
$$E=\bigcup_{i=1}^kE_i.$$
On the other hand
\begin{eqnarray*}K_Y^2&=&N_{Y|X}^2\\
&=&(-K_Z).E^3\\
&=&({\cal O}_Z(4)+p^*h).({\cal O}_Z(1)-p^*h)^3\\
&=&4c_1({\cal E})-12+1\\
&=&-7.
\end{eqnarray*}
Hence we arrive at
$$k=16.$$

Since $-K_Z=\phi^*{\cal O}_{{\bb P}^4}(5)-E'$ the characterization of $\phi|Y$ as a blowup shows, that
$X'$ is a quintic in ${\bb P}^4$ with $16$
double points on a linearly embedded ${\bb P}^2\subset{\bb P}^4$.


$(iv)$: In the last case $c_1({\cal E})=0$, obviously ${\cal E}=4{\cal O}$ and hence
$Z={\bb P}^1\times{\bb P}^3$. Therefore we may use Theorem \ref{3g<16} and get that $E$ consists
of $64$ curves
$C\cong{\bb P}^1$.
\end{proof}

\end{document}